\magnification=1200
\documentstyle{amsppt}
\raggedbottom
\hsize=16.35truecm
\vsize=24.5truecm
\NoRunningHeads
\NoBlackBoxes
\TagsOnRight
\output={\plainoutput}
\headline={\hfil}
\footline={\ifnum\count0=1\hfil\else\hss\folio\hss\fi}
\font\srm=cmr10 at 8pt
\font\sbf=cmb10 at 8pt
\font\ssl=cmsl10 at 8pt
\font\bbf=cmb10 at 12pt
\baselineskip=14pt
\rm

\redefine\!{\kern-.075em}
\def\qed{\hfill$\ssize\square$}
\predefine\accute{\'}
\redefine\'{\kern.05em{}}
\def\bbigskip{\vskip 24pt plus 8pt minus 8pt}
\def\list#1{\par\noindent\hangindent=1.5truecm\hangafter=-2
\hbox to 0pt{\hskip-\hangindent[#1]\hfill}}
\def\<{\langle}
\def\>{\rangle}
\def\Circ{\raise.05em\hbox{$\scriptstyle\,\circ\,$}}
\def\Ad{\operatorname{Ad}}
\def\id{\operatorname{id}}
\def\tr{\operatorname{Tr}}
\def\Ric{\operatorname{Ric}}
\def\grad{\operatorname{grad}}
\def\F{\Cal F}
\def\G{\Cal G}
\def\H{\Cal H}
\def\J{\Cal J}
\def\N{\Cal N}
\def\P{\Cal P}
\def\V{\Cal V}
\def\Z{\Cal Z}
\def\fs{\frak s}
\def\fu{\frak u}
\def\fC{\frak C}
\def\a{\alpha}
\def\b{\beta}
\def\e{\epsilon}
\def\s{\sigma}
\def\v{\varphi}
\def\nab#1{\nabla\kern-.2em\lower.8ex\hbox{$\ssize#1$}\'}
\def\diffop#1#2{#1\kern-.1em\lower.8ex\hbox{$\ssize#2$}\'}
\def\Rics#1{\hbox{\sl Ric\/}\kern-.05em\lower.6ex\hbox{$\ssize#1$}\,}
\def\dvs{d^{\'v}\kern-.15em\s}


\centerline{\bbf THE ENERGY OF UNIT VECTOR FIELDS ON THE 3-SPHERE}
\bigskip
\centerline{\bf by}
\bigskip
\centerline{\bbf A\. Higuchi, B\. S\. Kay \& C\. M\. Wood}
\bbigskip
{\sbf Abstract.} 
\srm 
The stability of the 3-dimensional Hopf vector field, as a harmonic
section of the unit tangent bundle, is viewed from a number of different
angles.  The spectrum of the vertical Jacobi operator is computed, and
compared with that of the Jacobi operator of the identity map on the
3-sphere. The variational behaviour of the 3-dimensional Hopf vector
field is compared and contrasted with that of the closely-related Hopf
map.  Finally, it is shown that the Hopf vector fields are the unique global
minima of the energy functional restricted to unit vector fields on the
3-sphere.
\bbigskip

\rm

\head 1. Introduction \endhead
A smooth unit vector field $\s$ on a compact Riemannian manifold $(M,g)$
with
Euler char\-acteristic zero may be regarded as a smooth mapping of
Riemannian
manifolds $\s\colon(M,g)\to(UM,h)$, where $UM$ is the unit tangent bundle
and
$h$ is the restriction of the Sasaki metric on the tangent bundle $TM$.
The energy of $\s$ may be defined accordingly.  Since metrics $h$ and $g$
are
horizontally isometric, and $\s$ is a section, it suffices to consider the
{\sl vertical energy functional:}
$$
E^v(\s)\,=\,\int_M|\'\dvs\'|^2\,dx
\tag1-1
$$
where $\dvs$ is the vertical component of the differential $d\s$.
(Here, `horizontal' and `vertical' refer to the complementary distributions
on $TM$ defined by the Levi-Civita connection).  One says that $\s$ is a
critical point of $E^v$, or a {\sl harmonic section\/} of $UM$, if $E^v$ is
stationary at $\s$ with respect to variations through unit vector fields.
The (non-linear) Euler-Lagrange equations for this variational problem are
\cite{17}:
$$
\nabla^*\nabla\s\,-\,|\'\nabla\s\'|^2\'\s\,=\,0,
\tag1-2
$$
where $\nabla^*\nabla$ is the {\sl trace\/} (or {\sl rough\/})
{\sl Laplacian:} 
$$
\nabla^*\nabla\s\,=\,-\tr\nabla^2\s.
$$
Further, one says that a harmonic section $\s$ is {\sl $E^v$-stable\/} if
the
second variation of $E^v$ at $\s$ with respect to unit vector fields is
non-negative.  The second variation of $E^v$ in this constrained sense may
be
regarded as a quadratic form $\H^v_\s$ (the {\sl vertical Hessian\/}) on the
space $\V_\s$ of appropriate variation fields; since $\s$ is allowed to
vary only through unit vector fields, $\V_\s$ is the space of smooth
vector fields on $M$ which are pointwise orthogonal to $\s$.  Associated to
$\H^v_\s$ is the {\sl vertical Jacobi operator\/} $\J^v_\s$:
$$
\H^v_\s(\a,\b)\,=\,\int_M\<\J^v_\s(\a),\b\>\,dx,
\qquad\text{for all $\a,\b\in\V_\s$.}
\tag1-3
$$
(Diamond brackets denote the relevant Riemannian metric, in this case $g$).
Explicit computation (see \cite{17}) shows that $\J^v_\s$ is the following
symmetric, elliptic linear second order partial differential operator on
$\V_\s$: 
$$
\J^v_\s(\a)\,=\,\nabla^*\nabla\a\,-\,|\'\nabla\s\'|^2\,\a\,
-\,2\<\nabla\s,\nabla\a\>\'\s,
\tag1-4
$$
Thus, $\J^v_\s$ may be viewed as a twisted version of the conventional
Jacobi
operator for a sphere-valued harmonic map \cite{14}; indeed,
since everything mentioned so far is in fact true for unit sections $\s$ of
any Riemannian vector bundle, this constitutes a complete generalization of
the theory of harmonic maps into spheres.  There is a unique extension of
$\J^v_\s$ to a self-adjoint linear operator on the Hilbert space of $L^2$
variation fields of $\s$.
\par
The {\sl canonical Hopf vector field\/} on $M=S^{2n+1}$ ($n=0,1,2,\dots$) is
defined: 
$$
\s(x)\,=\,ix,\qquad x\in\Bbb R^{2n+2}\cong\Bbb C^{n+1},\;|\'x\'|=1,
\tag1-5
$$
where $i=\sqrt{-1}$.  More generally, any unit vector field congruent to
$\s$ will be called a {\sl Hopf vector field.}  It is known that in all
dimensions $\s$ is a harmonic section, and if $n=2,3,\dots$ then $\s$ is
$E^v$-unstable (see \cite{17}).  This comes as no surprise, in view of
Y-L.\'Xin's instability theorem for harmonic mappings from spheres
\cite{19}.  More surprising is the fact that the $3$-dimensional Hopf vector
field is $E^v$-stable \cite{18}.  In this paper, we take a closer look at
$\J^v_\s$ when $\s$ is the 3-dimensional Hopf vector field, confirming its
non-negativity in a number of different ways.  In \S2 we make the
interesting observation that the non-negativity of $\J^v_\s$ is purely a
consequence of dimension and curvature, and the fact that $\s$ has geodesic
integral curves; the Lie group structure of $S^3$ is in fact only incidental
(although it was used in \cite{18}).  Our primary aim in \S2, however, is to
compute the (necessarily discrete) spectrum of $\J^v_\s$, along with the
eigenvalue multiplicities (Theorem 2.3).  It is interesting to compare these
spectral data with those of the Jacobi operator of the identity map
$\id$ on $S^3$; for, the Hessian $\H$ of the latter is related to $\H^v_\s$
by restriction:   
$$
\H^v_\s\,=\,\H(\V_\s,\V_\s).
$$
Recall that the space of variation fields for $\id$ is the entire Lie
algebra of vector fields on $S^3$; of course, since $\id$ is $E$-unstable,
$\H$ is indefinite.  Thus, we have a rather nice example illustrating how a
`simple' change of domain can dramatically alter the spectrum.  (See also
Remark 2.4).  In addition to this close relationship between the variational
aspects of $\s$ and $\id$, there is also a strong linkage between $\s$ and
another $E$-unstable harmonic map: the Hopf map  $\v\colon S^3\to S^2$.  At
first sight (see Proposition 3.1) this link seems so compelling that one is
tempted to infer that the variational properties of $\s$ must surely be the
same as those of $\v$.  In \S3 we tease out the subtleties of the
relationship, and show how $\s$ manages to avoid the instabilities of $\v$.
In addition, we compute the spectrum of the Jacobi operator of $\v$,
correcting errors of \cite{15} (see Remark 3.12).  Finally, in \S4 we
cap-off
the discussion of $\J^v_\s$ with the following much more powerful global
result:  
\medskip
\proclaim{Main Theorem}
The absolute minimum of $E^v$ over all unit vector fields on $S^3$ is
$2\pi^2$, which is achieved at, and only at, the Hopf vector fields.
\endproclaim 
\medskip\noindent
This thereom provides an affirmative resolution of a conjecture in
\cite{18}, and establishes a complete analogy between the behaviour of the
energy and volume functionals on the space of unit vector fields on
$S^3$ \cite{9}.  Our proof of the Main Theorem utilizes a rigidity theorem
for shear-free geodesic congruences on $S^3$, which follows
from the work of P.\,Baird and J.C.\,Wood on harmonic morphisms (combining
results of \cite{2} and \cite{3}).  However, we take the opportunity to
provide a direct proof of this rigidity theorem.


\head 2. The Spectrum of the Vertical Jacobi Operator \endhead
We refer to the following well-known diffeomorphism of $S^3$ with $SU(2)$ as
the {\sl Pauli correspondence:}
$$
(z_1,z_2)\,\longleftrightarrow\,
\pmatrix z_1&-\overline{z_2}\\z_2&\overline{z_1}\endpmatrix
\qquad z_1,z_2\in\Bbb C,\quad|\'z_1\'|^2+|\'z_2\'|^2=1,
$$
and the following basis of the Lie algebra $\fs\fu(2)$ as the
{\sl Pauli basis:}
$$
P_1=\pmatrix i&0\\0&-i\endpmatrix\qquad
P_2=\pmatrix 0&i\\i&0\endpmatrix\qquad
P_3=\pmatrix 0&-1\\1&0\endpmatrix
$$
The Pauli basis satisfies the commutation relations:
$$
[\'P_j,P_k\']\,=\,2\e_{jkl}\,P_l, \tag2-1
$$
where $\e_{jkl}$ is the totally antisymmetric symbol.
There exists a unique bi-invariant metric on $SU(2)$ rendering the Pauli
basis orthonormal, and with respect to this metric the Pauli
correspondence is an isometry.  Furthermore, with respect to the
orientation of $SU(2)$ determined by the Pauli basis, the Pauli
correspondence preserves the orientation of $S^3$ induced from the
standard orientation of $\Bbb R^4$ by the outward-pointing normal.
The Hopf vector fields on $S^3$ are precisely the Pauli-preimages of unit
vector fields on $SU(2)$ which are either left-invariant or
right-invariant.
\par
Let $\s_j$ denote the Pauli-preimage of the left-invariant vector field on
$SU(2)$ generated by $P_j$.  Then $(\s_1,\s_2,\s_3)$ is a
positively-oriented global orthonormal frame on $S^3$, with $\s_1=\s$, the
canonical Hopf vector field.  It follows from left-invariance (see Lemma
3.2) that 
$$
\nab{\s_j}\s_k\,=\,\frac12\'[\'\s_j,\s_k\']\,
=\,\e_{jkl}\,\s_l,\qquad\text{by (2-1).}
\tag2-2
$$
It follows from (2-2) that
$$
|\'\nabla\s\'|^2\,=\,2.
\tag2-3
$$
In order to compute the spectrum of $\J^v_\s$, we first rewrite this
operator using the triad basis $\{\s_j\}$, noting that any $\a\in\V_\s$ can
be written as $\,\a=f_2\s_2+f_3\s_3$, where $f_2$ and $f_3$ are smooth
real-valued functions on $S^3$.  The following fact will be useful.
\bigskip
\proclaim{2.1 Lemma}
If $X$ (resp\. $\lambda$) is any smooth vector field (resp\. function) on a
Riemannian manifold, then:
$$
\nabla^*\nabla(\lambda X)\,
=\,(\Delta\lambda)\'X\,+\,\lambda\'\nabla^*\nabla X\,
-\,2\'\nab{\grad\lambda}X,
$$
where $\Delta$ is the Laplace-Beltrami operator acting on functions.
\endproclaim
\bigskip 
\proclaim{2.2 Proposition}
Let
$$
\a\,=\,f_2\s_2\,+\,f_3 \s_3\,, \tag2-4
$$
where $f_2$ and $f_3$ are smooth real functions on $S^3$. Then
$$
\J^v_\s(\a)\,=\,\bigl(\Delta f_2\,+\,2\nab\s f_3\bigr)\s_2\,
+\,\bigl(\Delta f_3\,-\,2\nab\s f_2\bigr)\s_3\,.
$$
\endproclaim
\demo{Proof}
It follows from (1-4) and (2-3) that
$$
\J^v_\s(\a)\,=\,\nabla^*\nabla\a\,-\,2\a\,-\,2\<\nabla\s,\nabla\a\>\'\s.
\tag2-5
$$
We first compute $(\nabla^*\nabla-2)\a$.  By Lemma 2.1:
$$
\align
\nabla^*\nabla\a\,
&=\,(\Delta f_2)\s_2\,+\,(\Delta f_3)\s_3\,
+\,f_2\'\nabla^*\nabla\s_2\,+\,f_3\'\nabla^*\nabla\s_3 \\
&\quad-\,2\bigl(\nab{\grad f_2}\s_2\,
+\,\nab{\grad f_3}\s_3\bigr).
\endalign
$$
Since $\s_j$ is a Killing field, we have
$$
\nabla^*\nabla\s_j\,=\,\Ric(\s_j)\,=\,2\s_j.
$$
(This can also be deduced from equation (2-2)).  Furthermore, by (2-2):
$$
\align
\nab{\grad f_2}\s_2\,
&=\,\sum_{j=1}^3(\nab{\s_j}f_2)\,\nab{\s_j}\s_2\,
=\,\sum_{j,l}\e_{j2l}\,(\nab{\s_j}f_2)\,\s_l\,
=\,(\nab{\s_1}f_2)\,\s_3\,-\,(\nab{\s_3}f_2)\,\s_1\,, \\
\nab{\grad f_3}\s_3\,
&=\,\sum_{j,l}\e_{j3l}\,(\nab{\s_j}f_3)\,\s_l\,
=\,(\nab{\s_2}f_3)\,\s_1\,-\,(\nab{\s_1}f_3)\,\s_2\,. \\
\endalign
$$
Therefore:
$$
\align
(\nabla^*\nabla-2)\a\,
&=\,(\Delta f_2)\s_2\,+\,(\Delta f_3)\s_3\,
-\,2\bigl((\nab\s f_2)\s_3 - (\nab\s f_3)\s_2\bigr)\\
&\quad +\,2\'\bigl(\'\nab{\s_3}f_2 - \nab{\s_2}f_3\'\bigr)\'\s.
\tag2-6 
\endalign
$$
Next, we compute $\<\nabla\s,\nabla\a\>$.  By (2-2) again:
$$
\<\nabla\s,\nabla\a\>\,
=\,\sum_{j=1}^3 \<\nab{\s_j}\s,\nab{\s_j}\a\>\,
=\,-\<\s_3,\nab{\s_2}\a\>\,+\,\<\s_2,\nab{\s_3}\a\>.
$$
Hence
$$
\<\nabla\s,\nabla\a\>\,=\,\nab{\s_3}f_2 - \nab{\s_2}f_3\,.
$$
Therefore the term proportional to $\s$ in (2-6) is cancelled by
$-2\<\nabla\s,\nabla\a\>\s$ in (2-5) as it should.
Hence, the result follows. \qed
\enddemo
\bigskip
Now, let us define a map $\psi$ from the space of smooth vector
fields pointwise orthogonal to $\s$ to the space of smooth complex
functions on $S^3$ by
$$
\psi:\a\,\to\,f=f_2+if_3.
$$ 
The space of $\Bbb C$-valued functions on $S^3$ is regarded as a real
Hilbert space, using the inner product
$$
(f,g)\,=\,\operatorname{Re} \int_{S^3}f\bar{g}\,dx,
$$
where $\bar{f}$ is the complex conjugate of $f$.
Then the proposition implies that
$$
\J^v_\s\,=\,\psi^{-1}\circ\Lambda\circ\psi,
$$ 
where the map $\Lambda$ on the space of complex functions is defined by
$$
\Lambda(f)\,=\,\Delta f\,-\,2i\'\nab\s f.
$$
Hence, the spectrum of $\J^v_\s$ is identical to that of $\Lambda$.  Since
the operators $\Delta$ and $\nabla_{\s}$ commute, the complete set of
eigenfunctions of $\Lambda$ can be chosen to be simultaneous eigenfunctions
of $\Delta$ and$\nabla_{\s}$.  The eigenvalues of $\Delta$ are $n(n+2)$,
$n=0,1,2,\ldots$, and eigenfunctions corresponding to each $n$ form the
representation $(n/2,n/2)$ of $SU(2)\otimes SU(2)$, with multiplicity
$2(n+1)^2$.  (Here,  since the functions $f$ and $if$ are linearly
independent, the multiplicity is doubled compared to the usual case.)  It is
known that the eigenvalues of $i\'\nabla_{\s}$ on scalar functions are
$-n,-n+2,\ldots n-2,n$, each with multiplicity $2(n+1)$. Thus, we have
established the following theorem.
\bigskip
\proclaim{2.3 Theorem}
The eigenvalues of $\J^v_\s$, where $\s$ is a Hopf vector field, are
$n(n+2)+2k$, where $n$ is a non-negative integer and
$\,k=-n,-n+2,\dots,n-2,n$.  The multiplicity of each eigenvalue for given
$n$ and $k$ is $2(n+1)$.
\endproclaim
\bigskip
We note that two eigenvalues with different $n$ and $k$ are distinct.
In fact, the eigenvalue sequence is a splice of the two arithmetic
progressions $4k,4k+1$, $k=0,1,2,\dots$  Since the eigenvalues are all
non-negative, it follows that $\J^v_\s$ is non-negative.  The zero
eigenvalue
corresponds to $f=$ const. Thus, $\J^v_\s(\a)=0$ if and only if $\a$ is a
linear combination of $\s_2$ and $\s_3$; in particular, $\a$ is
left-invariant.
\bigskip
\subhead 2.4 Remark \endsubhead
The Jacobi operator $\J$ for the identity map on $S^3$ is \cite{14}:
$$
\J\,=\,\nabla^*\nabla-2\,=\,\Delta-4,
$$
by application of a Weitzenb\"ock formula, where $\Delta$ is the Hodge-de
Rham Laplacian acting on the entire space of 1-forms/vector fields.
The eigenvalues of $\Delta$ fall into two sequences:
$$
\lambda_k\,=\,(k+1)(k+3),\qquad
\tilde\lambda_k\,=\,(k+2)^2,
\qquad k=0,1,2,\dots
$$
with corresponding multiplicities:
$$
\mu_k\,=\,(k+2)^2,\qquad
\tilde\mu_k\,=\,2(k+1)(k+3),
\qquad k=0,1,2,\dots
$$
(The eigenvectors of $\lambda_k$ are exact $1$-forms, whereas those of
$\tilde\lambda_k$ are co-exact; see \cite{4} with corrections of \cite{12},
or \cite{7}).  Therefore the eigenvalues of $\J$ are
$$
k^2+4k-1,\qquad k^2+4k,\qquad k=0,1,2,\dots
$$
with multiplicities $\mu_k,\tilde\mu_k$ respectively.  Notice that these two
sequences are qualitatively quite different from those of Theorem 2.3.
\bigskip
The non-negativity of $\J^v_\s$ was established in a slightly
different manner in \cite{18} using the Bochner-Yano integral formula
\cite{5} on any compact manifold $(M,g)$:
$$
\int_M \left( |\nabla X|^2\,-\,\Ric(X,X)\right)\,dx\,
=\,\int_M \bigl(\tfrac12\'|\'\diffop LX g\'|^2\,-\,(\div X)^2\bigr)\,dx.
\tag2-7
$$
This was used to rewrite (1-3) on $S^3$ as
$$
\H^v_\s(\a,\a)\,
=\,\int_{S^3}\bigl(\tfrac12\'|\'\diffop L{\a}g\'|^2\,
-\,(\div\a)^2\bigr)\,dx,
\tag2-8
$$
if $\s$ is a Hopf vector field.  Then, the expansion (2-4) was used to show
the non-negativity of $\H^v_\s$, and to derive the condition on $\a$ for
$\H^v_\s(\a,\a) = 0$.  In fact, the non-negativity of the {\it integrand\/}
of (2-8) holds with a weaker condition on $\s$.  We conclude this section by
establishing this fact.  (A similar inequality will be
used in \S4 to prove the Main Theorem, mentioned in the Introduction).
\bigskip
\proclaim{2.5 Proposition (Linear Inequality)}
Let $\s$ be a unit vector field on an $n$-dimensional Riemannian manifold
$(M,g)$.  If the integral curves of $\s$ are geodesics, then for any
vector field $\a$ pointwise orthogonal to $\s$ we have:
$$
(n-1)\'|\'\diffop L{\a}g\'|^2\,\geqslant\,4(\div\a)^2.
$$
\endproclaim
\demo{Proof}
Recall first that
$$
\diffop LXg(Y,Z)\,
=\,\bigl\<\nab YX,Z\bigr\>\,+\,\bigl\<Y,\nab ZX\bigr\>
\tag2-9
$$
for any vector field $X$.  It follows that
$$
\tr\diffop LXg\,=\,2\div X.
\tag2-10
$$
Let $\a_1,\dots,\a_{n-1}$ be local vector fields such that
$(\a_1,\dots,\a_{n-1},\s)$ is an orthonormal moving frame on $M$.
It follows from (2-10) that:
$$
\align
(n-1)\'|\'\diffop LXg\'|^2\,&-\,4(\div X)^2\,
=\,(n-1)\'|\'\diffop LXg\'|^2\,-\,\bigl(\tr\diffop LXg\bigr)^2 \\
\vspace{1.5ex}
&=\,(n-2)\bigl(\diffop LXg(\s,\s)^2\,
+\,\sum_i\diffop LXg(\a_i,\a_i)^2\bigr)
+\,2(n-1)\sum_i\diffop LXg(\s,\a_i)^2 \\
\vspace{1.5ex}
&\quad+\,(n-1)\sum_{i\neq j}\diffop LXg(\a_i,\a_j)^2\,
-\,2\sum_i\diffop LXg(\s,\s)\'\diffop LXg(\a_i,\a_i)
\vspace{1ex}
&\qquad
-\,2\sum_{i<j}\diffop LXg(\a_i,\a_i)\'\diffop LXg(\a_j,\a_j).
\tag2-11
\endalign
$$
Since the integral curves of $\s$ are geodesics, it follows from (2-9) that:
$$
\diffop L{\a}g(\s,\s)\,=\,2\bigl\<\nab\s\a,\s\bigr\>\,
=-2\bigl\<\a,\nab\s\s\bigr\>\,=\,0.
\tag2-12
$$
Therefore if $X=\a$ then (2-11) collapses to:
$$
\align
(n-1)\'|\'\diffop LXg\'|^2\,-\,4(\div X)^2\,
&=\,2(n-1)\,\diffop LXg(\s,\a_i)^2 \,
+\,(n-1)\sum_{i\neq j}\diffop LXg(\a_i,\a_j)^2 \\
&\qquad+\sum_{i<j}\bigl(\diffop LXg(\a_i,\a_i)
-\diffop LXg(\a_j,\a_j)\bigr)^2
\tag2-13 \\
&\geqslant\,0.
\tag"$\ssize\square$"
\endalign
$$
\enddemo
\bigskip
Note that the non-negativity of the integrand in (2-8) follows from case
$n=3$ of this proposition.  Although there are many unit vector fields
on $S^3$ whose integral curves are geodesics \cite{8}, apart from the Hopf
vector fields we do not know of any which are harmonic sections of $US^3$.
In fact, this is a very interesting open question.


\head 3. Comparison with the Hopf Map\endhead
We first note an extremely simple relationship between the 3-dimensional
canonical Hopf vector field $\s$ and the Hopf map $\v\colon S^3\to S^2$.
Recall that on any Lie group $G$ there is the Maurer-Cartan form $\mu$,
with values in the Lie algebra:
$$
\mu(X)\,=\,g^{-1}.X,\qquad\forall X\in T_gG,\quad\forall g\in G.
$$
Let us denote by $\eta$ the right-invariant analogue:
$$
\eta(X)\,=\,X.g^{-1}
$$
Throughout \S3 we identify $S^3$ with $SU(2)$ via the Pauli correspondence,
as described at the beginning of \S2.  Then $\s$ is left-invariant;
therefore $\mu\Circ\s=$ const. (the Pauli matrix $P_1$).  On the other
hand: 
\bigskip
\proclaim{3.1 Proposition}
There exists an isometric identification of $S^2$ with the unit sphere in
$\fs\fu(2)$ such that $\eta\Circ\s=\v$, the Hopf map $S^3\to S^2$.
\endproclaim 
\demo{Proof}
If the Pauli-image of $(z_1,z_2)\in S^3$ is denoted by $\gamma$, and $e$
denotes the identity, then:
$$
\align
\eta\Circ\s(\gamma)\,&=\,\gamma.\s(e).\gamma^{-1}\,
=\,\pmatrix z_1&-\overline{z_2}\\z_2&\overline{z_1}\endpmatrix
\pmatrix i&0\\0&-i\endpmatrix
\pmatrix\overline{z_1}&\overline{z_2}\\-z_2&z_1\endpmatrix \\
&=\,\pmatrix i\'|\'z_1\'|^2-i\'|\'z_2\'|^2&2iz_1\'\overline{z_2}\\
2i\'\overline{z_1}\'z_2&i\'|\'z_2\'|^2-i\'|\'z_1\'|^2\endpmatrix \\
\vspace{1ex}
&=\,(|\'z_1\'|^2-|\'z_2\'|^2)P_1\,+\,2\Re(\overline z_1z_2)P_2\,
+\,2\Im(\overline z_1z_2)P_3.
\tag3-1
\endalign
$$
On the other hand, by the Hopf map we understand the composition of the
Hopf fibration $\pi\colon S^3\to\Bbb CP^1;(z_1,z_2)\mapsto(z_1\colon z_2)$
(homogeneous coordinates) with the standard isomorphism of $\Bbb CP^1$
with the unit sphere in $\Bbb R^3$.  This isomorphism admits the following
local description: on the open subset of $\Bbb CP^1$ where $z_1\neq0$, the
complex chart $(z_1\colon z_2)\mapsto z_2/z_1$ is followed by the inverse
of stereographic projection from the north pole onto the equatorial
plane.  Thus:
$$
\v(z_1,z_2)\,=\,(2\'\overline{z_1}\'z_2,\,|\'z_2\'|^2-|\'z_1\'|^2)\,
\in\,\Bbb C\times\Bbb R\cong\Bbb R^3.
$$
But by (3-1) this is precisely the map $S^3\to S^2$ obtained by taking the
coordinates of $\eta\Circ\s$ with respect to the basis $(P_2,P_3,-P_1)$
of $\fs\fu(2)$. \qed
\enddemo
\bigskip
\remark{Note}
The most natural identification of the unit sphere in $\fs\fu(2)$ with
$S^2$ (viz\. the restriction of the linear isomorphism which sends the
Pauli basis $(P_1,P_2,P_3)$ to the standard basis of $\Bbb R^3$)
is not the identification asserted by Proposition 3.1; rather, the two
differ by a rotary reflection of $S^2$.
\endremark
\bigskip
At first glance, the fact that $\s$ differs from $\v$ `only' by right
translation suggests a complete equivalence between the variational
theories of $E^v(\s)$ and $E(\v)$.  Of course, in the light of \S2 we know
this is not the case, since $\v$ is $E$-unstable.  A closer examination of
the energy functionals (see Proposition 3.4), and Jacobi operators (see
Proposition 3.7) reveals a more complicated  scenario, which resolves
the issue.  Recall first the following characterizations of covariant
differentiation in Lie groups \cite{10}.
\bigskip   
\proclaim{3.2 Lemma}
For any Lie group with bi-invariant metric, the Levi-Civita
connection may be characterized in either of the following two ways:
\smallskip\noindent
{\rm(a)}\quad
$\mu\bigl(\nab XY\bigr)\,=\,d(\mu Y)(X)\,+\,\tfrac12\'[\'\mu X,\mu Y\']$
\smallskip\noindent
{\rm(b)}\quad
$\eta\bigl(\nab XY\bigr)\,
=\,d(\eta Y)(X)\,-\,\tfrac12\'[\'\eta X,\eta Y\']$.
\endproclaim
\bigskip
\proclaim{3.3 Lemma}
If $Y$ is a vector field on any Lie group with bi-invariant metric, then:
$$
|\'d(\eta Y)(X)\'|^2\,=\,2\'|\nab XY\'|^2\,-\,|\'d(\mu Y)(X)\'|^2\,
+\,\tfrac12\'|\'[\'\mu X,\mu Y\']\'|^2.
$$
In particular, on $S^3\cong SU(2)$ we have:
$$
|\'d(\eta Y)\'|^2\,=\,4\'|\'Y\'|^2\,+\,2\'|\'\nabla Y\'|^2\,
-\,|\'d(\mu Y)\'|^2.
$$
\endproclaim
\demo{Proof}
It follows from Lemma 3.2\,(b) and the bi-invariance of the metric that:
$$
\align
|\'d(\eta Y)(X)\'|^2\,
&=\,|\'\nab XY\'|^2\,
+\,\bigl\<[\'\eta X,\eta Y\'],\,\eta\bigl(\nab XY\bigr)\bigr\>\,
+\,\tfrac14\'|\'[\'\eta X,\eta Y\']\'|^2 \\
&=\,|\'\nab XY\'|^2\,
+\,\bigl\<[\'\mu X,\mu Y\'],\,\mu\bigl(\nab XY\bigr)\bigr\>\,
+\,\tfrac14\'|\'[\'\mu X,\mu Y\']\'|^2.
\endalign
$$
Now use the algebraic fact that if $a,b,c$ are vectors in any inner
product space, with $a=b+c$, then
$$
2\<a,c\>\,=\,|\'a\'|^2\,-\,|\'b\'|^2\'+\'|\'c\'|^2.
$$
Lemma 3.2\,(a) permits the choice
$$
a\,=\,\mu\bigl(\nab XY\bigr),\quad
b\,=\,d(\mu Y)(X),\quad
c\,=\,\tfrac12\'[\'\mu X,\mu Y\']
$$
and the general identity follows. This identity implies:
$$
|\'d(\eta Y)\'|^2\,=\,2\'|\'\nabla Y\'|^2\,-\,|\'d(\mu Y)\'|^2\,
+\,\tfrac12\sum_j|\'[\'\mu E_j,\mu Y\']\'|^2
$$
where $\{E_j\}$ is any orthonormal tangent frame.  On $SU(2)$, let
$\{E_j\}$ be the global left-invariant orthonormal frame $\{\s_j\}$, and
write $Y=Y^k\s_k$ (summation convention).  By the commutation relations
(2-1):  
$$
\sum_j|\'[\'\mu E_j,\mu Y\']\'|^2\,
=\,\sum_{j,k}|\'[\'P_j,\,Y^kP_k\']\'|^2\,
=\,\sum_{j,k,l}4\'|\'\e_{jkl}\,Y^kP_l\'|^2\,
=\,8\sum_k(Y^k)^2\,=\,8\'|\'Y\'|^2.
\tag"$\ssize\square$"
$$
\enddemo
\bigskip
\subhead Remark\endsubhead
It follows from Proposition 3.1 and Lemma 3.3 (taking $Y=\s$) that:
$$
|\'d\v\'|^2\,=\,4\,+\,2\'|\'\nabla\s\'|^2\,=\,8,
$$
and we recover the familiar fact that the Hopf map $\v\colon S^3\to S^2$
is an `eigenmap' with eigenvalue $8$ \cite{6}.
\bigskip
\proclaim{3.4 Proposition}
Let $\s_t$ be any variation of the Hopf vector field on $S^3$ through unit
vector fields, and let $\v_t=\eta\Circ\s_t$ be the corresponding variation
of the Hopf map (cf\. Proposition 3.1).  Then:
$$
E(\v_t)\,=\,4\pi^2\,+\,2\'E^v\!(\s_t)\,-\,E(\mu\Circ\s_t).
$$
\endproclaim
\demo{Proof}
Integrate Lemma 3.3 with $Y=\s_t$, recalling that
$\,|\'\dvs_t\'|^2=|\'\nabla\s_t\'|^2$ (see \cite{17}), and $S^3$ has
volume $2\pi^2$.  \qed
\enddemo
\bigskip
It follows from Proposition 3.4 that if $\s_t$ is $E^v$-decreasing, then
$\v_t$ is $E$-decreasing.  However, the converse is not necessarily true.
Indeed, since Proposition 3.4 tells us that $E^v$ is essentially
the average of the energies of the right- and left-translates, it is
conceivable that any $E$-decreasing variation of the right-translate $\v$
of $\s$ is compensated by an $E$-increasing variation of the
left-translate.  Our aim is to show that this is indeed the case.
\par
Before proceeding further, we make some simple observations.
Let $\<\s\>$ denote the line subbundle of $TS^3$ generated by $\s$.  Then
$\<\s\>=\ker d\v$.  Therefore $\V_\s$ is precisely the space of smooth
{\sl $\v$-horizontal\/} vector fields on $S^3$.  This terminology will be
very useful.  For example, the covariant derivative of $\s$ may be
written as follows:
$$
\nab X\s\,=\,\cases
0,\quad\text{if $X$ is $\v$-vertical} \\
iX,\quad\text{if $X$ is $\v$-horizontal.} \endcases
\tag3-2
$$
\par
The space of $L^2$ variation fields for $\v$ admits the following
$L^2$-orthogonal decomposition:
$$
\N_\v\oplus\Z_\v\oplus\P_\v
$$
where $\N_\v$ (resp\. $\P_\v\'$) is the direct sum of the negative (resp\.
positive) eigenspaces of the Jacobi operator $\J_\v$ of $\v$ (see
\cite{14}), and $\Z_\v$ is the kernel of $\J_\v$.  General elliptic theory
guarantees that $\N_\v$ and $\Z_\v$ are finite-dimensional, and that all
their elements are smooth sections (of the pullback bundle $\v^*TS^2$).
The dimensions of $\N_\v$ and $\Z_\v$ were computed in \cite{15} (modulo a
few errors; see Remark 3.12), from which the following facts may be
deduced. 
\smallskip\noindent
(A)\quad
$\N_\v$ is 4-dimensional; it comprises variation fields of the form
$d\v(\Gamma)$ where $\Gamma$ is a conformal gradient field on $S^3$.
(By a conformal gradient field on a sphere $S^n$ we mean the spherical
gradient of the restriction of a linear functional on the ambient
Euclidean space $\Bbb R^{n+1}$).  Let $\G_\s$ denote the subspace of
$\V_\s$ comprising the $\v$-horizontal components of conformal gradient
fields on $S^3$ (no conformal gradient field is $\v$-horizontal).  Then
$\N_\v=d\v(\G_\s)$.
\smallskip\noindent
(B)\quad
$\Z_\v$ is 8-dimensional; it is generated by variation fields of the form
$d\v(Z)$ where $Z$ is an infinitesimal isometry of $S^3$, and of the form
$C(\v)$ where $C$ is a conformal vector field on $S^2$.  The Lie algebra
$\frak I$ of infinitesimal isometries of $S^3$ is 6-dimensional; it is
generated by the vector fields on $SU(2)$ which are either left- or
right-invariant.  However, the fibres of $\v$ are invariant under the flow
of the left-invariant vector field $\s$, so $d\v(\frak I)$ is actually
5-dimensional.  The Lie algebra $\fC$ of conformal vector fields on
$S^2$ is also 6-dimensional.  We may write $C(\v)=d\v(\tilde C)$ where
$\tilde C$ is the $\v$-horizontal lift of $C$.  If $C$ is an infinitesimal
isometry of $S^2$ then $\tilde C$ is an infinitesimal isometry of $S^3$
(in fact, a right-invariant vector field), so in this case
$C(\v)\in d\v(\frak I)$.  However, if $C$ is a conformal gradient field of
$S^2$ then the variation fields $C(\v)$ constitute a 3-dimensional subspace
of $\fC(\v)$ which is complementary to $d\v(\frak I)$ in $\Z_\v$.
(Note that $\tilde C$ is not a conformal field on $S^3$, unless $C$ is an
infinitesimal isometry).  If $\tilde\fC$ denotes the $6$-dimensional
subspace of $\V_\s$ comprising the $\v$-horizontal lifts of elements of
$\fC$, and an $8$-dimensional subspace $\F_\s$ of $\V_\s$ is defined
by: 
$$
\F_\s\,=\,\tilde\fC\oplus\Bbb R\s_2\oplus\Bbb R\s_3,
\tag3-3
$$
then it follows that $\Z_\v=d\v(\F_\s)$.
\smallskip 
The space of $L^2$ variation fields for $\s$ also decomposes as an
$L^2$-orthogonal direct sum:
$$
\N_\s\oplus\Z_\s\oplus\P_\s
$$
where $\N_\s$ (resp\. $\P_\s\'$) is the direct sum of the negative
(resp\. positive) eigenspaces of $\J^v_\s$, and $\Z_\s$ is the kernel
of $\J^v_\s$.  We would like to see how the subspaces $\F_\s$ and $\G_\s$
corresponding to low spectral frequencies of $\J_\v$ relate to the
eigenspaces of $\J^v_\s$ (see Proposition 3.8); when used in conjunction
with our energy formula (Proposition 3.4), this will enable us to deduce
that $\N_\s$ is trivial.   We will derive this relationship by comparing
the Jacobi operators $\J_\v$ and $\J^v_\s$ (see Proposition 3.7).
\par
Recall that the Jacobi operator $\J_f$ for an arbitrary harmonic mapping
$f$ of Riemannian manifolds is given by \cite{14}:
$$
\J_f(w)\,=\,\nabla^*\nabla w\,-\,\Rics{f}(w),
$$
where $w$ is any variation field for $f$, and
$$
\Rics{f}(w)\,=\,\tr R(w,df)df.
$$
If the domain of $f$ is a sphere $S^{m+1}$, and the variation field is of
the form $w=df(X)$ where $X$ is a vector field on $S^{m+1}$, then
application of a Weitzenb\"ock formula yields \cite{19}:
$$
\J_f(w)\,=\,df(\nabla^*\nabla X-m\'X)\,
-\,2\sum_j\nabla df\bigl(E_j,\nab{E_j}X\bigr),
\tag3-4
$$
where $\{E_j\}$ is any local orthonormal tangent frame in $S^{m+1}$.
The term $\nabla df$ is sometimes referred to as the {\sl second
fundamental form\/} of $f$; when $f=\v$ it can be computed using standard
results on Riemannian submersions.
\bigskip
\proclaim{3.5 Lemma}
Let $\a,\b$ (resp\. $V,W$) be $\v$-horizontal (resp\. $\v$-vertical)
tangent vectors of $S^3$.  Then:
\smallskip\noindent
{\rm(a)}\quad
$\nabla d\v(\a,\b)=0$;
\smallskip\noindent
{\rm(b)}\quad
$\nabla d\v(V,W)=0$;
\smallskip\noindent
{\rm(c)}\quad
$\bigl\<\nabla d\v(\a,V),\,d\v(\b)\bigr\>\,
=\,2\'\<V,\'[\'\a,\b\']\>$.
\endproclaim
\demo{Proof}
The map $\tfrac12\v$ is a Riemannian submersion.  (By $k\v$ for any
$k\in\Bbb R^+$ we mean the mapping onto the $2$-sphere of radius $k$
obtained by scalar multiplication in ambient $\Bbb R^3$).  Then, (a) is an
identity for all Riemannian submersions, (b) follows from the fact
that $\v$ has totally geodesic fibres, and (c) is a rescaling of the
identity
$$
\bigl\<\nabla d\pi(\a,V),\,d\pi(\b)\bigr\>\,
=\,\tfrac12\'\<V,\'[\'\a,\b\']\>
$$
for any Riemannian submersion $\pi$.  (Standard results on the second
fundamental form of a Riemannian submersion may be found in
\cite{1,\,Ch.\,9}, \cite{11}, \cite{13}, \cite{16}). \qed
\enddemo
\bigskip
Let $J$ denote the standard complex structure on the unit sphere $S^2$,
characterized by the condition that $\v$ is `horizontally homothetically
holomorphic': 
$$
d\v(i\a)\,=\,2J\'d\v(\a),
\tag3-5
$$
for all $\v$-horizontal vectors $\a$.
\bigskip
\proclaim{3.6 Lemma}
For all vector fields $X,Y$ on $S^3$ we have:
$$
\nabla d\v(X,Y)\,=\,-2J\'d\v\bigl(\<\s,X\>Y+\<\s,Y\>X\bigr).
$$
\endproclaim
\demo{Proof}
Note first that if $\a,\b$ are $\v$-horizontal then:
$$
\align
\<\s,[\'\a,\b\']\>\,&=\,\bigl\<\s,\nab\a\b-\nab\b\a\bigr\>\,
=\,-\bigl\<\nab\a\s,\b\bigr\>+\bigl\<\nab\b\s,\a\bigr\> \\
&=\,-2\<i\a,\b\>,\qquad\text{by (3-2).}
\endalign
$$
If $\{\b_j:j=1,2\}$ is any $\v$-horizontal orthonormal frame,
then $\{\tfrac12\'d\v(\b_j)\}$ is an orthonormal frame of $S^2$,
and it follows from Lemma 3.5 and (3-5) that:
$$
\align
\nabla d\v(\s,\a)\,
&=\,\tfrac14\'\<\nabla d\v(\s,\a),\,d\v(\b_j)\>\'d\v(\b_j)\,
=\,\tfrac12\'\<\s,\'[\'\a,\b_j\']\>\'d\v(\b_j) \\
&=\,-\<i\a,\b_j\>\'d\v(\b_j)\,=\,-d\v(i\a)\,=\,-2J\'d\v(\a).
\tag3-6
\endalign
$$
Now suppose that $\a,\b$ are the $\v$-horizontal components of $X,Y$
respectively:
$$
X\,=\,\a\,+\,\<X,\s\>\s,\qquad Y\,=\,\b\,+\,\<Y,\s\>\s.
$$
\par\pagebreak\noindent
It follows from Lemma 3.5 that:
$$
\allowdisplaybreaks
\align
\nabla d\v(X,Y)\,&=\,\<X,\s\>\'\nabla d\v(\s,\b)\,
+\,\<Y,\s\>\'\nabla d\v(\a,\s) \\
&=\,-2\'\<\s,X\>\'J\'d\v(\b)\,
-\,2\'\<\s,Y\>\'J\'d\v(\a),\qquad\text{by (3-6)} \\
&=\,-2J\'d\v\bigl(\<\s,X\>Y+\<\s,Y\>X\bigr).
\tag"$\ssize\square$"
\endalign
$$
\enddemo
\bigskip
\proclaim{3.7 Proposition}
Let $\a$ be any $\v$-horizontal vector field on $S^3$, and let
$w=d\v(\a)$.  Then:
$$
\J_\v(w)\,=\,d\v\bigl(\J^v_\s(\a)\bigr)\,-\,2w\,
+\,4\'J\'d\v\bigl(\nab\s\a\bigr).
$$
\endproclaim
\demo{Proof}
It follows from (3-4) and (2-5) that
$$
\J_\v(w)\,=\,d\v\bigl(\J^v_\s(\a)\bigr)\,
-\,2\sum_j\nabla d\v\bigl(E_j,\nab{E_j}\a\bigr).
$$
By Lemma 3.6 (applying the summation convention):
$$
\align
-2\,\nabla d\v\bigl(E_j,\nab{E_j}\a\bigr)\,
&=\,4\'J\'d\v\bigl(\<E_j,\s\>\'\nab{E_j}\a\,
+\,\bigl\<\nab{E_j}\a,\s\bigr\>\'E_j\bigr) \\
&=\,4\'J\'d\v\bigl(\nab\s\a\,
-\,\bigl\<\a,\nab{E_j}\s\bigr\>E_j\bigr) \\
&=\,4\'J\'d\v\bigl(i\a\,+\,\nab\s\a\bigr),
\qquad\text{by (3-2),} \\
&=\,-2\'d\v(\a)\,+\,4\'J\'d\v\bigl(\nab\s\a\bigr),
\qquad\text{by (3-5).}
\tag"$\ssize\square$"
\endalign
$$
\enddemo
\bigskip
\proclaim{3.8 Proposition}
\flushpar
{\rm(a)}\quad
$\F_\s$ and $\G_\s$ are complex subspaces of $\V_\s$.
\smallskip\noindent
{\rm(b)}\quad
$\F_\s\subset\Z_\s\oplus\P_\s$.
\smallskip\noindent
{\rm(c)}\quad
$\G_\s\subset\P_\s$.
\endproclaim
\remark{Note}
Since the $\v$-horizontal distribution is a complex vector bundle, $\V_\s$
is a complex vector space.
\endremark
\demo{Proof}\;(a)\quad
It was shown in \cite{18} that $\,i\G_\s=\G_\s$.  On the other hand, $\fC$
may be viewed as the Lie algebra of holomorphic vector fields on the
Riemann sphere, and is therefore a complex Lie algebra.  It follows from
(3-5) that $i\tilde\fC=\tilde\fC$.  Furthermore $i\s_2=-\s_3$, so the
subspace $\Bbb R\s_2\oplus\Bbb R\s_3$ is also complex.  It therefore
follows from (3-3) that $i\F_\s=\F_\s$.
\smallskip\noindent
(b)\quad
If $\a$ is a $\v$-basic vector field on $S^3$ (ie\. a $\v$-horizontal
field which projects to $S^2$), then $[\'\s,\a\']$ is $\v$-adapted
to the zero vector field on $S^2$, and hence $\v$-vertical.
Therefore:  
$$
d\v\bigl(\nab\s\a\bigr)\,=\,d\v\bigl(\nab\a\s\bigr)\,=\,d\v(i\a).
$$
If in addition $d\v(\a)$ is a Jacobi field for $\v$ then Proposition 3.7
reads: 
$$
d\v\bigl(\J^v_\s(\a)\bigr)\,
=\,2\'d\v(\a)\,-\,4\'J\'d\v(i\a)\,=\,4\'d\v(\a),
\qquad\text{by (3-2),}
$$
which implies $\,\J^v_\s(\a)=4\a$.  In particular, this shows that
$\tilde\fC\subset\P_\s$.  On the other hand, if $\a$ is any
left-invariant $\v$-horizontal vector field then $\,\J^v_\s(\a)=0$,
since any such vector field is Hopf.   (See also
Proposition 2.2, with $f_2$ and $f_3$ constant).
\smallskip\noindent
(c)\quad
It was shown in \cite{18} that every element of $\G_\s$ is an eigenvector
of $\J^v_\s$ with eigenvalue $1$.  (This can also be seen from
Proposition 3.7, since the variation field $w=d\v(\a)$ is a $-1$
eigenvector for $\J_\v$\'---\'see \cite{19}\'---\'and $\nab\s\a$ is
$\v$-vertical.  In fact, since $\G_\s$ is $4$-dimensional, it follows from
Theorem 2.3 that $\G_\s$ is the entire eigenspace with eigenvalue $1$).
\qed 
\enddemo
\bigskip
A variation $\s_t$ of $\s$ through unit vector fields produces a variation
$\v_t=\eta\Circ\s_t$ of the Hopf map, by Proposition 3.1.  However, as
noted in (A) and (B) above, the most natural variations of $\v$ with
variation fields in $\N_\v$ or $\Z_\v$ are of the form $\v\Circ\psi_t$,
where $\{\psi_t\}$ is the flow of a ($\v$-horizontal) vector field on
$S^3$.  In order to play off Proposition 3.8 against Proposition 3.4 we
need to relate the variation fields of these two types of variation (see
Proposition 3.11). 
\bigskip
\proclaim{3.9 Lemma}
If $\a$ is any $\v$-horizontal vector field on $S^3$ then
$\,d\v(\a)=2\eta(i\a)$.
\endproclaim
\demo{Proof}
It follows from Lemma 3.2\,(b), with $Y=\s$, and Proposition 3.1 that
$$
d\v(\a)\,=\,\eta\bigl(\nab\a\s\bigr)\,
+\,\tfrac12\'\Ad(g)\'[\'\mu\Circ\a,P_1\'].
$$
Write $\,\a=f_2\s_2+f_3\s_3$, so that $\,\mu\Circ\a=f_2P_2+f_3P_3$.  The
commutation relations (2-1) yield:
$$
[\'\mu\Circ\a,P_1\']\,=\,[\,f_2P_2+f_3P_3,\,P_1\']\,
=\,2f_3P_2-2f_2P_3\,=\,2\'\mu(i\a).
$$
Using (3-2), it follows that:
$$
d\v(\a)\,=\,\eta(i\a)\,+\,\Ad(g)\'\mu(i\a)\,=\,2\'\eta(i\a).
\tag"$\ssize\square$"
$$
\enddemo
\bigskip
\proclaim{3.10 Lemma}
If $K$ is the connection map for any linear connection on a Lie group $G$,
then
$$
d\eta(V)\,=\,\eta(KV),\qquad\text{for all vertical $V\in T(TG)$.}
$$
\endproclaim
\demo{Proof}
Let $g\in G$, and $X,Y\in T_gG$.
Suppose $V$ is the `vertical lift' of $Y$ at $X$:
$$
V\,=\,\left.\frac{d}{dt}\right|_{t=0}(X+tY).
$$
Then $KV=Y$, and therefore $\eta(KV)=Y.g^{-1}$.  On the other hand:
$$
d\eta(V)\,=\,\left.\frac{d}{dt}\right|_{t=0}\bigl((X+tY).g^{-1}\bigr)\,
=\,\left.\frac{d}{dt}\right|_{t=0}(X.g^{-1}+tY.g^{-1})\,
=\,Y.g^{-1}.
\tag"$\ssize\square$"
$$
\enddemo
\bigskip
\proclaim{3.11 Proposition}
Let $\s_t$ be a variation of the canonical Hopf vector field on $S^3$
through unit  vector fields, with variation field $\a$.  Define two
variations of the Hopf map $\v$ as follows:
$$
\v_t\,=\,\eta\Circ\s_t\qquad\text{and}\qquad
\Phi_t\,=\,\v\Circ\psi_t,
$$
where $\{\psi_t:t\in\Bbb R\}$ is the flow of the vector field
$X=-i\a$ on $S^3$.  If $w$ (resp\. $W$) is the variation field of $\v_t$
(resp\. $\Phi_t$) then $W=2w$.
\endproclaim
\demo{Proof}
Let $V$ denote the vertical lift of $\a$ into the tangent bundle of
$TS^3$.  Thus, for each $x\in S^3$, $V\bigl(\s(x)\bigr)$ is the element of
$T_{\s(x)}TS^3$ tangent to the curve $\s_t(x)$ in $T_xS^3$ at $t=0$.
Furthermore $\a=KV$, by a characteristic property of connection
maps.  By Lemma 3.10, the variation field for $\v_t$ is:
$$
w\,=\,d\eta(V)\,=\,\eta(KV)\,=\,\eta(\a).
$$
On the other hand, by Lemma 3.9:
$$
W\,=\,d\v(X)\,=\,-d\v(i\a)\,=\,2\'\eta(\a).
\tag"$\ssize\square$"
$$
\enddemo
\bigskip
To show that $\N_\s$ is trivial, we argue by contradiction.
Suppose that $\s_t$ is a variation of $\s$ through unit vector fields,
with variation field $\a$ a non-trivial negative eigenvector of
$\J^v_\s$.  It follows from parts (b) and (c) of Proposition 3.8 that $\a$
is $L^2$-orthogonal to $\F_\s\oplus\G_\s$.  Hence, by part (a) of
Proposition 3.8, $X=-i\a$ is also $L^2$-orthogonal to $\F_\s\oplus\G_\s$.
Therefore $d\v(X)$ is $L^2$-orthogonal to $\Z_\v\oplus\N_\v$, since
elements of $\F_\s$ and $\G_{\s}$ are $\v$-horizontal, and $\v$ is
horizontally homothetic.  Now $d\v(X)$ is the variation field $W$ for the
variation $\Phi_t=\v\Circ\psi_t$ of $\v$, where $\{\psi_t\}$ is the flow of
$X$.  Since $W\in\P_\v$, it follows from Proposition 3.11 that
$w\in\P_\v$ also, where $w$ is the variation field for
$\v_t=\eta\Circ\s_t$.  Therefore, by the second variation formula for
harmonic maps, $\v_t$ is $E$-increasing for small $|\'t\'|$.  However
$\s_t$ is $E^v$-decreasing in a neighbourhood of $t=0$, and it therefore
follows from Proposition 3.4 that $\v_t$ is $E$-decreasing.
\bigskip
\subhead 3.12 Remark\endsubhead
The energy formula (Proposition 3.4) may be used to compute the spectrum
of the Jacobi operator $\J_\v$ of the Hopf map, thereby correcting
errors of \cite{15}.  Define $\,E^1(\s_t)=E(\eta\Circ\s_t)$ and
$\,E^2(\s_t)=E(\mu\Circ\s_t)$, and let $\J^1_\s$ and $\J^2_\s$ denote the
corresponding Jacobi operators at the Hopf vector field $\s$; note from
Proposition 3.1 that $\J^1_\s$ is conjugate to $\J_\v$.  It follows from
Proposition 3.4 that
$$
\J^1_\s\,=\,2\J^v_\s\,-\,\J^2_\s.
$$
Write $\a\in\V_\s$ as $\,\a=f_2\s_2+f_3\s_3$, and define
$\,\psi\colon\a\mapsto f=f_2+if_3$, as in \S2.  Then $f$ is a variation
field for the constant map $\mu\Circ\s$, whose Jacobi operator is the
Laplace-Beltrami operator:
$$
\psi\Circ\J^2_{\s}\Circ\psi^{-1}(f)\,=\,\Delta f.
$$
It therefore follows from Proposition 2.2 that
$$
\psi\Circ\J^1_\s\Circ\psi^{-1}(f)\,
=\,2(\Delta f\,-\,2i\,\nab\s f)\,-\,\Delta f\,
=\,\Delta f\,-\,4i\,\nab\s f.
$$
By arguing as in the proof of Theorem 2.3, it follows that the eigenvalues
of
$\J_\v$ are $\,n(n+2)+4k$ where $n$ is a non-negative integer and
$\,k=-n,-n+2,\dots,n-2,n$, with multiplicity $2(n+1)$.  In particular, the
only negative eigenvalue is $n=1,k=-1$; thus $\N_\v$ is $4$-dimensional.
Furthermore, the only possibilities for eigenvalue zero are $n=0,k=0$ and
$n=2,k=-2$.  The multiplicity of the former is $2$, and that of the latter
is $2(2+1)=6$; thus $\Z_\v$ is indeed $8$-dimensional.


\head 4. Proof of the Main Theorem\endhead
First, recall that a vector field $X$ on a Riemannian manifold $(M,g)$
is said to define a {\sl conformal foliation,} or {\sl shear-free
congruence,} if 
$$
\diffop LXg(A,B)\,=\,\lambda\,g(A,B),
\tag4-1
$$
for all vector fields $A,B$ pointwise orthogonal to $X$, where
$\lambda\colon M\to\Bbb R$ is a smooth function.  (It follows from
(2-10) that $\,(n-1)\lambda=2\div\s$, where $n$ is the dimension of $M$).
If in addition the integral curves of $X$ are (possibly reparametrized)
geodesics, then $X$ defines a {\sl conformal geodesic foliation,} or
{\sl shear-free geodesic congruence.}  Such vector fields appear naturally
in the following non-linear version of Proposition 2.5.
\bigskip
\proclaim{4.1 Proposition (Non-Linear Inequality)}
Let $\s$ be a unit vector field on an $n$-dimensional Riemannian
manifold.  Then
$$
(n-1)\'|\'\diffop L{\s}g\'|^2\,\geqslant\,4(\div\s)^2,
$$
with equality if and only if $\s$ is a shear-free geodesic congruence.
\endproclaim
\demo{Proof}
The argument used to derive this inequality is almost identical to that of
Proposition 2.5; however, in place of identity (2-12), equation (2-9)
may be used to derive:
$$
\diffop L{\s}g(\s,\s)\,=\,2\bigl\<\nab\s\s,\s\bigr\>\,
=\,\s.|\'\s\'|^2\,=\,0.
$$
(Note that it is not necessary to assume that the integral curves of $\s$
are geodesics).  Thus equation (2-13) holds with $X=\s$.
\par
Since $\s$ has constant length, it follows from (2-9) that
$$
\diffop L{\s}g(\s,\a_i)\,=\,\bigl\<\nab\s\s,\a_i\bigr\>,
$$
so the vanishing of these Lie derivatives implies $\nab\s\s=0$.  On the
other hand, the vanishing of the remaining Lie derivatives in (2-13)
implies
$$
\diffop L{\s}g(\a,\b)\,=\,\lambda\,g(\a,\b),
$$
for all $\a,\b$ orthogonal to $\s$, where $\lambda\colon M\to\Bbb R$ is a
smooth function.  Thus, equality  occurs only when $\s$ defines a
conformal geodesic foliation. \qed
\enddemo
\bigskip
Our main theorem can now be proved by using the Bochner-Yano integral
formula (2-7) and the following consequence of \cite{2} and \cite{3}:
\bigskip 
\proclaim{Rigidity Theorem}
A unit vector field $\s$ is a shear-free geodesic congruence on
$S^3$ if and only if $\s$ is a Hopf vector field.
\endproclaim
\bigskip
Note that the ``if'' part is obvious.
We shall give a direct proof of this theorem later for
completeness.  Let us first prove our Main Theorem using this result.
\bigskip
\proclaim{Main Theorem}
The absolute minimum of $E^v$ over all unit vector fields on $S^3$ is
$2\pi^2$, which is achieved at, and only at, the Hopf vector fields.
\endproclaim
\demo{Proof}
By the Bochner-Yano integral formula (2-7), we find on $S^3$ that
$$
E^v(\s)\,=\,\frac14\int_{S^3}
\left( |\'\diffop L{\s}g\'|^2\,-\,2(\div\s)^2\right)\,dx\,+\,2\pi^2.
$$
We have used $\,\Ric(\s,\s)=2$ and the fact that the volume of
$S^3$ is $2\pi^2$.
By Proposition 4.1 with $n=3$, we have $E^v\!(\s)\geqslant 2\pi^2$ and
that the minimum is achieved if and only if $\s$ is a shear-free geodesic
congruence.  Then, the theorem follows immediately from the Rigidity
Theorem. \qed 
\enddemo
\bigskip
We conclude this paper by proving the ``only if'' part of
the Rigidity Theorem as promised.  We first prove the following lemma.
\bigskip
\proclaim{4.2 Lemma}
Suppose that $\s$ is a shear-free congruence on a Riemannian
$3$-manifold.  Let $\a$, $\b$ and $\s$ form a local orthonormal
frame, and define $Z=\a+i\b$.  Then
$$
\nab Z\s\,=\,\Phi\,Z,
\tag4-2
$$
where $\Phi\colon M\to\Bbb C$ is a smooth complex function.
\endproclaim
\demo{Proof}
By the definition (4-1) of shear-free congruences we obtain:
$$
\<\a,\nab\a\s\>\,=\,\<\b,\nab\b\s\>\,=\,\lambda/2\,=\,\mu,
\quad\text{say,}
$$
and
$$
\<\a,\nab\b\s\>\,+\,\<\b,\nab\a\s\>\,=\,0.
$$
By defining $\nu=\<\a,\nab\b\s\>$ we find:
$$
\nab\a\s\,=\,\mu\'\a-\nu\'\b\,,\qquad
\nab\b\s\,=\,\mu\'\b+\nu\'\a.
$$
Then, we obtain the desired formula by letting $\,\Phi=\mu+i\nu$. \qed
\enddemo
\bigskip
A couple of differential equations can be derived for $\Phi$
on a manifold of constant curvature.
\bigskip
\proclaim{4.3 Proposition}
Let $\s$ be a shear-free geodesic congruence on a $3$-manifold
of constant curvature $c$.  Then
$$
\align
\s.\Phi\,&=\,-c-\Phi^2, \tag4-3\\
\bar Z.\Phi\,&=\,0, \tag4-4
\endalign
$$
where $\Phi$ and $Z$ are defined in Lemma 4.2.
\endproclaim
\demo{Proof}
>From (4-2) we find
$$
\nab\s(\nab Z\s)\,=\,(\s.\Phi)Z\,+\,\Phi\,\nab\s Z.
\tag4-5
$$
On the other hand we have
$$
\align
\nab\s(\nab Z\s)\,
&=\,R(\s,Z)\s\,+\,\nab Z(\nab\s\s)\,
+\,\nab{\nab\s Z}\s\,-\,\nab{\nab Z\s}\s \\
&=\,c(\<Z,\s\>\s-Z)\,+\,\nab Z(\nab\s\s)\,
+\,\nab{\nab\s Z}\s\,-\,\nab{\nab Z\s}\s.
\tag4-6
\endalign
$$
Note that for any vector $X$,
$$
X\,=\,\frac12\left(\<\bar Z,X\>\'Z\,
+\,\<Z,X\>\'\bar Z\right)\,+\,\<\s,X\>\'\s.
$$
Using this formula, $\nab\s\s= 0$ and equation (4-2), we find from (4-6):
$$
\nab\s(\nab Z\s)\,=\,\Phi\,\nab\s Z\,-\,(c+\Phi^2)Z.
$$
By comparing this with (4-5) we obtain (4-3).

Next we have
$$
\align
\bigl(\nab Z\nab{\bar Z}\,-\,\nab{\bar Z}\nab Z\bigr)\s\,
&=\,\nab Z(\bar\Phi\bar Z)\,-\,\nab{\bar Z}(\Phi Z) \\
&=\,(Z.\bar\Phi)\'\bar Z\,+\,\bar\Phi\,\nab Z\bar Z\,
-\,(\bar Z.\Phi)\'Z\,-\,\Phi\,\nab{\bar Z}Z.
\tag4-7
\endalign
$$
On the other hand
$$
\align
(\nab Z\nab{\bar Z}\,-\,\nab{\bar Z}\nab Z)\'\s\,
&=\,R(Z,\bar Z)\s\,+\,\nab{\nab Z\bar Z}\s\,
-\,\nab{\nab{\bar Z}Z}\s \\
&=\,\bar\Phi\,\nab Z\bar Z\,-\,\Phi\,\nab{\bar Z}Z.
\endalign
$$
By comparing this formula with (4-7) we find
$$
(Z.\bar\Phi)\'\bar Z\,-\,(\bar Z.\Phi)\'Z\,=\,0.
$$
By taking the inner product with $\bar Z/2$ we obtain (4-4). \qed
\enddemo
\bigskip
This proposition allows us to compute the Laplacian of $\Phi$, as shown in
the  following lemma.
\bigskip
\proclaim{4.4 Lemma}
Suppose $\s$ is a shear-free geodesic congruence on a  $3$-manifold of
constant curvature.  Then the function $\Phi$ is harmonic.
\endproclaim
\demo{Proof}
Note first that
$$
-\Delta\Phi\,=\,\nab\s\nab\s\Phi\,+\,\nab Z\nab{\bar Z}\Phi\,
-\,\nab{\nab Z\bar Z}\Phi.
$$
Then, using Proposition 4.3 and
$$
\align
\nab Z\bar Z\,
&=\,\bigl\<\s,\nab Z\bar Z\bigr\>\'\s\,
+\,\tfrac12\bigl(\bigl\<\bar Z,\nab Z\bar Z\bigr\>\'Z\,
+\,\bigl\<Z,\nab Z\bar Z\bigr\>\'\bar Z\bigr) \\
&=\,-\bigl\<\nab Z\s,\bar Z\bigr\>\'\s\,
+\,\tfrac12\bigl\<Z,\nab Z\bar Z\bigr\>\'\bar Z \\
&=\,-2\'\Phi\'\s\,+\,\tfrac12\bigl\<Z,\nab Z\bar Z\bigr\>\'\bar Z,
\endalign
$$
we find $\Delta\Phi=0$. \qed
\enddemo
\bigskip
The Rigidity Theorem follows from this lemma.
\bigskip
\demo{Proof of the Rigidity Theorem}
Since $S^3$ is compact and $\Phi$ is harmonic, $\Phi$ is constant.
Then, from (4-3) we find $\Phi=\pm i$ on $S^3$ with $c=1$.  Thus,
$$
\nab Z\s\,=\,\pm iZ,  \tag4-8
$$
which implies that $\s$ is a Hopf vector field.  This fact can be seen as
follows.  Equation (4-8) can be written as
$$
\nab\a\s\,=\,\mp\b,\qquad \nab\b\s\,=\,\pm\a.
$$
These equations, together with $\nab\s\s=0$ imply
$$
\diffop L{\s}g(X,Y)\,=\,\bigl\<X,\nab Y\s\bigr\>\,
+\,\bigl\<Y,\nab X\s\bigr\>\,=\,0,
$$
for any vector fields $X$ and $Y$. Thus, $\s$ is a Killing vector field
of unit length, hence a Hopf vector field. \qed
\enddemo


\centerline{\bf References.}
\srm
\medskip
\list{1}A.\,Besse,
`Einstein Manifolds',
Ergebnisse der Mat\. und ihrer Grenzgebiete {\sbf 10},
Springer Verlag, 1987.
\smallskip
\list{2}P.\,Baird \& J.C.\,Wood,
{\ssl Bernstein theorems for harmonic morphisms from $\Bbb R^3$ and
$S^3$,} 
Math\. Ann\. {\sbf 280} (1988), 579--603.
\smallskip
\list{3}P.\,Baird \& J.C.\,Wood,
{\ssl Harmonic morphisms and conformal foliations by geodesics of
three-dimensional space forms,}
J\. Austral\. Math\. Soc\. {\sbf 51} (1991), 118--153.
\smallskip
\list{4}R.\,Camporesi \& A.\,Higuchi,
The Plancherel measure for $p$-forms in real hyperbolic spaces,
J\. Geom\. Physics {\sbf 15} (1994), 57--94.
\smallskip
\list{5}K.\,Yano \& S.\,Bochner,
`Curvature and Betti Numbers',
Ann\. of Math\. Studies {\sbf 32}, Princeton Univ\. Press, 1953.
\smallskip
\list{6}J.\,Eells \&  A.\,Ratto,
`Harmonic Maps and Minimal Immersions with Symmetries',
Annals of Math\. Studies {\sbf 130}, Princeton Univ\. Press, 1993.
\smallskip
\list{7}S.\,Gallot \& D.\,Meyer,
{\ssl Operateur de courbure et Laplacien des formes differentielles
d'une vari\accute et\accute e Riemannienne,}
J\. Math\. Pure Appl\. {\sbf 54} (1975), 259--284.
\smallskip
\list{8}H.\,Gluck \& F.W.\,Warner,
{\ssl Great circle fibrations of the three-sphere,}
Duke Math\. J\. {\sbf 50} (1983), 107--132.
\smallskip
\list{9}H.\,Gluck \& W.\,Ziller,
{\ssl On the volume of a unit vector field on the three-sphere,}
Comment\. Math\. Helv\. {\sbf 61} (1986), 177--192.
\smallskip
\list{10}S.\,Helgason,
`Differential Geometry, Lie Groups, and Symmetric Spaces',
Academic Press, 1978.
\smallskip
\list{11}R.\,Hermann,
{\ssl A sufficient condition that a mapping of Riemannian manifolds be a
fibre bundle,}  
Proc\. Amer\. Math\. Soc\. {\sbf 11} (1960), 236--242.
\smallskip
\list{12}A.\,Ikeda \& Y.\,Taniguchi,
{\ssl Spectra and eigenforms of the Laplacian on $S^n$ and
$P^n(\Bbb C)$,} Osaka J\. Math\. {\sbf 15} (1978), 515--546.
\smallskip
\list{13}B.\,O'Neill,
{\ssl The fundamental equations of a submersion,}
Michigan Math\. J\. {\sbf 13} (1966), 459--469.
\smallskip
\list{14}R.\,T.\,Smith,
{\ssl The second variation formula for harmonic mappings,}
Proc\. Amer\. Math\. Soc\. {\sbf 47} (1975), 229--236.
\smallskip 
\list{15}H.\,Urakawa,
{\ssl Stability of harmonic maps and eigenvalues of the Laplacian,}
Trans\. Amer\. Math\. Soc\. {\sbf 301} (1987), 557--589.
\smallskip
\list{16}J.\,Vilms,
{\ssl Totally geodesic maps,}
J\. Diff\. Geometry {\sbf 4} (1970), 73--79.
\smallskip 
\list{17}C.\,M.\,Wood,
{\ssl On the energy of a unit vector field,}
Geom\. Dedicata {\sbf 64} (1997), 319--330.
\smallskip
\list{18}C.\,M.\,Wood,
{\ssl The energy of Hopf vector fields,}
Manuscripta Math\., to appear.
\smallskip
\list{19}Y-L.\,Xin,
{\ssl Some results on stable harmonic maps,}
Duke Math\. J\. {\sbf 47} (1980), 609--613.
\bbigskip
{\hfill\rm Department of Mathematics, University of York,
Heslington, York Y01 5DD, U.K.}
\medskip
{\hfill\ssl E-mails:\quad\srm
ah28\@york.ac.uk, bsk2\@york.ac.uk, cmw4\@york.ac.uk}
\medskip
{\hfill\ssl 1991 Mathematics Subject Classification:\quad\srm
53C20 (53C15, 58E15, 58E20).}

\enddocument
\end